\documentclass{ifacconfN}

\usepackage{graphicx}      
\usepackage{natbib}        

\usepackage{amsmath}
\usepackage{amsthm}
\usepackage{amssymb}

\newtheorem{dfn}{Definition}
\newtheorem{rmk}{Remark}
\newtheorem{thmm}{Theorem}
\newtheorem{prp}[thmm]{Proposition}

\begin{document}
\begin{frontmatter}

\title{Port-Hamiltonian systems on discrete manifolds}

\author[a]{Marko \v{S}e\v{s}lija}
\author[a]{Jacquelien M.A. Scherpen}
\author[b]{Arjan van der Schaft}
\address[a]{Department of Discrete Technology and Production Automation, Faculty of Mathematics and Natural Sciences, University of Groningen, Nijenborgh 4, 9747 AG Groningen, The Netherlands, e-mail:~\{M.Seslija,\;J.M.A.Scherpen\}@rug.nl}
\address[b]{Johann Bernoulli Institute for Mathematics and Computer Science, University of Groningen, Nijenborgh 9, 9747 AG Groningen, The Netherlands, e-mail:~A.J.van.der.Schaft@rug.nl}

\begin{abstract}
This paper offers a geometric framework for modeling port-Hamiltonian systems on discrete manifolds. The simplicial Dirac structure, capturing the topological laws of the system, is defined in terms of primal and dual cochains related by the coboundary operators. This finite-dimensional Dirac structure, as discrete analogue of the canonical Stokes-Dirac structure, allows for the formulation of finite-dimensional port-Hamiltonian systems that emulate the behaviour of the open distributed-parameter systems with Hamiltonian dynamics.
\end{abstract}


\end{frontmatter}

\section{Introduction}
A large class of open distributed-parameter Hamiltonian systems can be defined with respect to the Stokes-Dirac structure \cite{vdSM02}.
This infinite-dimensional Dirac structure provides a theoretical account that permits the inclusion of varying boundary variables in the boundary problem for partial differential equations. From an interconnection and control viewpoint, such a treatment of boundary conditions is essential for the incorporation of energy exchange through the boundary, since in many applications the interconnection with the environment takes place precisely through the boundary. For numerical integration, simulation and control synthesis, it is of paramount interest to have finite approximations that can be interconnected to one another.

Most of the numerical algorithms emanating from the field of numerical analysis and scientific computing, however, fail to capture the intrinsic system structures and properties, such as symplecticity, conservation of momenta and energy, as well as differential gauge symmetry. Furthermore, some important results, including the Stokes theorem, fail to apply numerically and thus lead to spurious results. 

Recently in \cite{SeslijaCDC}, we have suggested a discrete exterior geometry approach to structure-preserving discretization of distributed-parameter port-Hamiltonian systems.  The spatial domain in the continuous theory represented by a finite-dimensional smooth manifold is replaced by a homological manifold-like simplicial complex and its circumcentric dual. The smooth differential forms, in discrete setting, are mirrored by cochains on the primal and dual complexes, while the discrete exterior derivative is defined to be the coboundary operator. A discrete analogue of the Stokes-Dirac structure is a so-called simplicial Dirac structure defined over a space of primal and dual discrete differential forms.

In this paper, we address the issue of matrix representations of simplicial Dirac structures by representing cochains by their coefficient vectors. In this manner, all linear operator from the continuous world can be represented by matrices, including the Hodge star, the coboundary and a trace operator. First, we recall the definition of the Stokes-Dirac structure and port-Hamiltonian systems. In the third section, we define some essential concepts from discrete exterior calculus as developed in \cite{Desbrun1,Desbrun, Hirani}. In order to allow the inclusion of nonzero boundary conditions on the dual cell complex, in \cite{SeslijaCDC} we have adapted a definition of the dual boundary operator that leads to a discrete analogue of the integration by parts formula, which is a crucial ingredient in establishing simplicial Dirac structures on a primal simplicial complex and its circumcentric dual. Finally, we demonstrate how these simplicial Dirac structures relate to spatially discretized wave equation on a bounded domain and the telegraph equations on a segment.
 
\section{The Stokes-Dirac structure and port-Hamiltonian systems}
The Stokes-Dirac structure is an infinite-dimensional Dirac structure that provides a foundation for port-Hamiltonian formulation of a class of distributed-parameter systems with boundary energy flow \cite{vdSM02}. 

Throughout this paper, let $M$ be an oriented $n$-dimensional smooth manifold with a smooth $(n-1)$-dimensional boundary $\partial M$ endowed with the induced orientation, representing the space of spatial variables. By $\Omega^k(M)$, $k=0,1,\ldots,n$, denote the space of exterior $k$-forms on $M$, and by $\Omega^k(\partial M)$, $k=0,1,\ldots,n-1$, the space of $k$-forms on $\partial M$. 

For any pair $p,q$ of positive integers satisfying $p+q=n+1$, define the flow and effort linear spaces by
\begin{equation*}
\begin{split}
\mathcal{F}_{p,q}=\,&\Omega^p(M)\times \Omega^q(M)\times \Omega^{n-p}(\partial M)\,\\
\mathcal{E}_{p,q}=\,&\Omega^{n-p}(M)\times \Omega^{n-q}(M)\times \Omega^{n-q}(\partial M)\,.
\end{split}
\end{equation*}
The bilinear form on the product space $\mathcal{F}_{p,q}\times\mathcal{E}_{p,q}$ is
\begin{equation}\label{eq-aj9}
\begin{split}
\langle\!\langle &(\underbrace{f_p^1, f_q^1, f_b^1}_{\in \mathcal{F}_{p,q}},\underbrace{e_p^1, e_q^1, e_b^1}_{\in \mathcal{E}_{p,q}}), (f_p^2, f_q^2, f_b^2, e_p^2, e_q^2, e_b^2)\rangle\!\rangle \\
& =\int_M e_p^1 \wedge f_p^2+e_q^1\wedge f_q^2+ e_p^2\wedge f_p^1+e_q^2\wedge f_q^1\\
&~~~ + \int_{\partial M} e_b^1 \wedge f_b^2 + e_b^2 \wedge f_b^1.
 \end{split}
\end{equation}

\vspace{0.1cm}
\begin{thmm}
Given linear spaces $\mathcal{F}_{p,q}$ and $\mathcal{E}_{p,q}$, and the bilinear form $\langle\!\langle, \rangle\!\rangle$, define the following linear subspace $\mathcal{D}$ of $\mathcal{F}_{p,q}\times \mathcal{E}_{p,q}$\vspace{-0.1cm}
\begin{equation}\label{eq-7Cont}
\begin{split}
\mathcal{D}=\big\{& (f_p,f_q,f_b,e_p,e_q,e_b)\in \mathcal{F}_{p,q}\times \mathcal{E}_{p,q}\big |\\
& \left(\begin{array}{c}f_p \\f_q\end{array}\right)=\left(\begin{array}{cc}0 & (-1)^{pq+1} {\mathrm{d}}\\ {\mathrm{d}} & 0\end{array}\right)\left(\begin{array}{c}e_p \\e_q\end{array}\right)\,,\\
& \left(\begin{array}{c}f_b \\e_b\end{array}\right)=\left(\begin{array}{cc}1&0\\0 & -(-1)^{n-q}\end{array}\right)\left(\begin{array}{c}e_p|_{\partial K} \\e_q|_{\partial K}\end{array}\right)
\big \}\,,
 \end{split}
\end{equation}
where $\mathrm{d}$ is the exterior derivative and $|_{\partial M}$ stands for a trace on the boundary $\partial M$. 
Then $\mathcal{D}=\mathcal{D}^\perp$, that is, $\mathcal{D}$ is a Dirac structure.
\end{thmm}


Consider a Hamiltonian density $\mathcal{H}:\Omega^p(M) \times \Omega^q(M)\rightarrow \Omega^n(M)$ resulting with the Hamiltonian $H=\int_M \mathcal{H}\in \mathbb{R}$. Setting the flows $f_p=-\frac{\partial \alpha_p}{\partial t}$, $f_q=-\frac{\partial \alpha_q}{\partial t}$ and the efforts $e_p=\delta_p H$, $e_q=\delta_q H$, where $(\delta_pH,\delta_qH)\in \Omega^{n-p}(M)\times \Omega^{n-q}(M)$ are the variational derivatives of $H$ at $(\alpha_p,\alpha_q)$, the distributed-parameter port-Hamiltonian system is defined by the relation
$$\left(  -\frac{\partial \alpha_p}{\partial t},  -\frac{\partial \alpha_q}{\partial t},f_b, \delta_p H,\delta_qH, e_b\right)\in \mathcal{D}\,,~~t\in \mathbb{R}\,.$$
Since $\frac{\textmd{d} H}{\textmd{d} t}=\int_{\partial M}e_b\wedge f_b$, the system is lossless.

\section{Basics of discrete exterior calculus}     

In the discrete setting, the smooth manifold $M$ is replaced by an oriented manifold-like simplicial complex. An $n$-dimensional simplicial complex $K$ is a simplicial triangulation of an $n$-dimensional polytope $|K|$ with an $(n-1)$-dimensional boundary. Familiar examples of such a discrete manifold are meshes of triangles embedded in $\mathbb{R}^3$ and tetrahedra obtained by tetrahedrization of a $3$-dimensional manifold.


\subsection{Chains and cochains}
The discrete analogue of a smooth $k$-form on the manifold $M$ is a $k$-cochain on the simplicial complex $K$. A $k$-chain is a formal sum of $k$-simplices of $K$ such that its value on a simplex changes sign when the simplex orientation is reversed. The free Abelian group generated by a basis consisting of oriented $k$-simplices with real-valued coefficients is $C_k(K;\mathbb{R})$. The space $C_k(K;\mathbb{R})$ is a vector space with dimension equal to the number of $k$-simplices in $K$, which is denoted by $N_k$. The space of $k$-cochains is the vector space dual of $C_k(K;\mathbb{R})$ denoted by $C^k(K;\mathbb{R})$ or $\Omega_d^k(K)$, as a reminder that this is the space of discrete $k$-forms. The space $\Omega_d^k(K)$ is the space of real-valued linear functionals on the vector space $C_k(K;\mathbb{R})$.



The discrete exterior derivative $\mathbf{d}^k:C^k(K)\rightarrow C^{k+1}(K)$ is defined by duality to the boundary operator $\partial_{k+1}:C_{k+1}(K;\mathbb{Z})\rightarrow C_{k}(K;\mathbb{Z})$, with respect to the natural pairing between discrete forms and chains. For a discrete form $\alpha\in \Omega_d^k(K) $ and a chain $c_{k+1}\in C_{k+1}(K;\mathbb{Z})$ we define $\mathbf{d}^k$ by\vspace{-0.2cm}
\begin{equation*}
\langle \mathbf{d}^k \alpha,c_{k+1}\rangle=\langle \alpha,(\mathbf{d}^k)^\textsc{t} c_{k+1}\rangle=\langle \alpha,\partial_{k+1}c_{k+1}\rangle\,,
\end{equation*}
where the boundary operator $\partial_{k+1}$ is the incidence matrix from the space of $k+1$-simplices to the space of $k$-simplices and is represented by a sparse $N_{k+1}\times N_{k}$ matrix containing only $0$ or $\pm1$ elements \cite{Desbrun2}. The important property of the boundary operator is $\partial_{k}\circ \partial_{k+1}=0$. The exterior derivative as the coboundary operator also satisfies $\mathbf{d}^{k+1}\circ\mathbf{d}^k =0$, what is a discrete analogue of the vector calculus identities $\mathrm{curl}\circ \mathrm{grad}=0$ and $\mathrm{div} \circ \mathrm{curl}=0$.


\subsection{Dual cell complex}
An essential ingredient of discrete exterior calculus is the dual complex of a manifold-like simplicial complex. Given a simplicial well-centered complex $K$, we define its interior dual cell complex $\star_\mathrm{i} K$ (block complex in terminology of algebraic topology \cite{Munkres}) as a circumcentric dual restricted to $|K|$. An important property of the the Voronoi duality is that primal and dual cells are orthogonal to each other. The boundary dual cell complex $\star_\mathrm{b} K$ is a dual to $\partial K$. The dual cell complex $\star K$ is defined as $\star K=\star_\mathrm{i} K \times \star_\mathrm{b} K$. A dual mesh $\star_\mathrm{i} K$ is a dual to $K$ in sense of a graph dual, and the dual of the boundary is equal to the boundary of the dual, that is $\partial (\star K)=\star (\partial K)=\star_\mathrm{b} K$. Because of duality, there is a one-to-one correspondence between $k$-simplices of $K$ and interior $(n-k)$-cells of $\star K$. Likewise, to every $k$-simplex on $\partial K$ there is a uniquely associated $(n-1-k)$-cell on $\partial(\star K)$.  Fig.~1 illustrates the duality on a flat $2$-dimensional simplicial complex.


\begin{figure}
\centering
    \includegraphics[width=8.4cm]{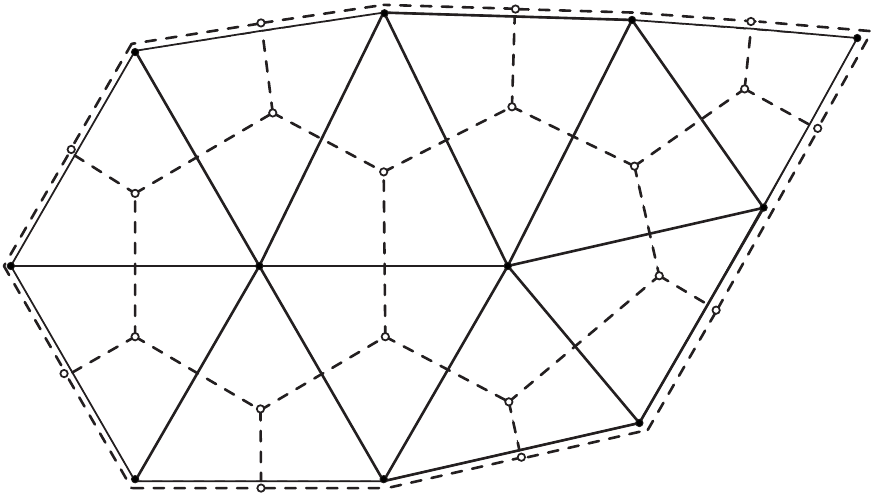}\\
  \caption{A $2$-dimensional simplicial complex $K$ and its circumcentric dual cell complex $\star K$ indicated by dashed lines. The boundary of $\star K$ is the dual of the boundary of $K$.}\label{fig:RDcomp1}
\end{figure}

In order to properly account for the behaviours on the boundary, we need to adapt the definition of the boundary dual operator as presented in \cite{Hirani, Desbrun}. In \cite{SeslijaCDC} we have proposed the following definition.
\vspace{0.2cm}
\begin{dfn}
The dual boundary operator $\partial_k: C_k(\star_\mathrm{i} K;\mathbb{Z})\\\rightarrow C_{k-1}(\star K;\mathbb{Z})$ is a homomorphism defined by its action on a dual cell $\hat{\sigma}^k=\star_\mathrm{i} \sigma^{n-k}=\star_\mathrm{i}[v_0,\ldots,v_{n-k}]$,
\begin{equation*}
\begin{split}
\partial_k \hat{\sigma}^k=\partial \star_{\mathrm{i},k} [v_0,\ldots,v_{n-k}]&=\partial_{\mathrm{i},k} \star_\mathrm{i} [v_0,\ldots,v_{n-k}]\\
&~~+\partial_{\mathrm{b},k} \star_\mathrm{i} [v_0,\ldots,v_{n-k}]\,,
\end{split}
\end{equation*}
where
\begin{equation*}
\begin{split}
\partial_{\mathrm{i},k} \star_\mathrm{i} [v_0,\ldots,v_{n-k}]&=\sum_{\sigma^{n-k+1}\succ \sigma^{n-k}}\star_\mathrm{i} (s_{\sigma^{n-k+1}}\sigma^{n-k+1})\\
\partial_{\mathrm{b},k} \star_\mathrm{i} [v_0,\ldots,v_{n-k}]&=\star_\mathrm{b} \left(s_{\sigma^{n-k}}\sigma^{n-k} \right)\,.
\end{split}
\end{equation*}

\end{dfn}

The boundary of the dual cell complex as defined in \cite{Hirani} is equal to $\partial_\mathrm{i}$. The dual boundary operator $\partial_\mathrm{b}$ extends the definition from \cite{Hirani} in such a manner that the boundary of the extended dual cell complex $\star K$ is the geometric boundary. The dual exterior derivatives $\mathbf{d}_\mathrm{i}^{k-1}:C^{k-1}(\star_\mathrm{i} K)\rightarrow C^{k}(\star_\mathrm{i} K)$ and $\mathbf{d}_\mathrm{b}^{k-1}:C^{k-1}(\star_\mathrm{b} K)\rightarrow C^{k}(\star_\mathrm{i} K)$ are defined by duality to the dual boundary operators $\partial_{\mathrm{i},k}$ and $\partial_{\mathrm{b},k}$, respectively.

\subsection{Discrete wedge and Hodge operator}
There exists a natural pairing, via the so-called primal-dual wedge product, between a primal $k$-cochain and a dual $(n-k)$-cochain. The resulting discrete form is the volume form. Let $\alpha^k\in \Omega_d^k(K)$ and $\hat{\beta}^{n-k}\in \Omega_d^{n-k}(\star_\mathrm{i} K)$. We define the discrete primal-dual wedge product $\wedge: \Omega_d^k(K)\times \Omega_d^{n-k}(\star_\mathrm{i} K) \rightarrow \Omega_d^{n}(V_k(K))$ by $
\langle \alpha^k\wedge \hat{\beta}^{n-k},V_{\sigma^{k}}\rangle
=\langle \alpha^k,\sigma^k \rangle \langle \hat{\beta}^{n-k}, \star_\mathrm{i}\sigma^k\rangle 
=(-1)^{k(n-k)}\langle \hat{\beta}^{n-k}\wedge \alpha^k,V_{\sigma^{k}}\rangle$, where $V_{\sigma^k}$ is the $n$-dimensional support volume obtained by taking the convex hull of the simplex $\sigma^k$ and its dual $\star \sigma^k$.

The proposed definition of the dual boundary operator ensures the validity of the evaluation by parts relation that parallels the integration by parts formula for smooth differential forms. 

\vspace{0.1cm}
\begin{prp}\label{prop:prim-dual}
Let $K$ be an oriented well-centered simplicial complex. Given a primal $(k-1)$-form $\alpha$ and an internal dual $(n-k)$-discrete form $\hat{\beta}_\mathrm{i} \in \Omega_d^{n-k}(\star_\mathrm{i}K)$ and a dual boundary form $\hat{\beta}_\mathrm{b}\in \Omega_d^{n-k}(\star_\mathrm{b}K) $, then
\begin{equation*}
\begin{split}
\langle \mathbf{d}^{k+1} \alpha \wedge \hat{\beta}_\mathrm{i} ,K\rangle\!&+\!(-1)^{k-1}\!\langle \alpha \!\wedge \!( \mathbf{d}_\mathrm{i}^{n-k} \hat{\beta}_\mathrm{i} \!+\!\mathbf{d}_\mathrm{b}^{n-k} \hat{\beta}_\mathrm{b} ),K\rangle \\
&~~~~~~~~~~~~~~~~~~~=\langle  \mathbf{tr}^{k-1}\alpha \wedge \hat{\beta}_\mathrm{b} ,\partial K\rangle\,,
\end{split}
\end{equation*}
where $\mathbf{tr}^{k-1}$ is the trace operator that isolates the components of the primal $(k-1)$-cochain associated with the boundary $\partial K$.\end{prp}

From definition of the dual boundary operator follows that $\mathbf{d}_\mathrm{i}^{n-k}=(-1)^k (\mathbf{d}^{k-1})^\textsc{t}$ and $\mathbf{d}_\mathrm{b}^{n-k}=(-1)^{k-1} (\mathbf{tr}^{k-1})^\textsc{t}$ (confer to \cite{SeslijaArXiv}).

The support volumes of a simplex and its dual cell are the same, which suggests that there is a natural identification between primal $k$-cochains and dual $(n-k)$-cochains. In the exterior calculus for smooth manifolds, the Hodge star, denoted $*_k$, is an isomorphism between the space of $k$-forms and $(n-k)$-forms. The discrete Hodge star is a map $*_k : \Omega_d^k(K)\rightarrow \Omega_d^{n-k}(\star_\mathrm{i} K)$ defined by its value over simplices and their duals. The Hodge star $*_k$ is a diagonal $N_k\times N_k$ matrix with the entry corresponding to a simplex $\sigma^k$ being $|\sigma^k|/|\star \sigma^k|$.

Another possibility for the construction of the Hodge operator us to use Whitney forms. The Whitney map is an interpolation scheme for cochains. It maps discrete forms to square integrable forms that are piecewise smooth on each simplex. The Whitney maps are built from barycentric coordinate functions and the resulting matrix is sparse but in general not diagonal \cite{Bosavit,Hiptmair}.




\section{Simplicial Dirac structures}\label{Sec4}
In this section, we introduce Dirac structures defined in terms of primal and duals cochains on the underlying discrete manifold. We call these Dirac structures \emph{simplicial Dirac structures}.

In the discrete setting, the smooth manifold $M$ is replaced by an $n$-dimensional well-centered oriented manifold-like simplicial complex $K$. The flow and the effort spaces will be the spaces of complementary primal and dual forms. The elements of these two spaces are paired via the discrete primal-dual wedge product. Let
\begin{equation*}
\mathcal{F}_{p,q}^d=\Omega_d^p(\star_\mathrm{i} K)\times \Omega_d^q( K)\times \Omega_d^{n-p}(\partial (K))\,
\end{equation*}
and 
\begin{equation*}
\mathcal{E}_{p,q}^d=\Omega_d^{n-p}( K)\times \Omega_d^{n-q}(\star_\mathrm{i} K)\times \Omega_d^{n-q}(\partial (\star K))\,.
\end{equation*}

The primal-dual wedge product ensures a bijective relation between the primal and dual forms, between the flows and efforts. A natural discrete mirror of the bilinear form (\ref{eq-aj9}) is a symmetric pairing on the product space $\mathcal{F}_{p,q}^d\times \mathcal{E}_{p,q}^d$ defined by
\begin{equation}\label{eq:bild}
\begin{split}
\langle\!\langle (&\underbrace{\hat{f}_p^1,{f}_q^1,{f}_b^1}_{\in \mathcal{F}_{p,q}^d},\underbrace{{e}_p^1,\hat{e}_q^1,\hat{e}_b^1}_{\in \mathcal{E}_{p,q}^d}), (\hat{f}_p^2,{f}_q^2,{f}_b^2,{e}_p^2,\hat{e}_q^2,\hat{e}_b^2)\rangle\!\rangle_d\\
&= \langle {e}_p^1\wedge \hat{f}_p^2+\hat{e}_q^1\wedge {f}_q^2+ {e}_p^2\wedge \hat{f}_p^1+\hat{e}_q^2\wedge {f}_q^1,K \rangle\\
&~\;~+\langle \hat{e}_b^1\wedge {f}_b^2+ \hat{e}_b^2\wedge f_b^1,\partial K \rangle
  \,.
  \end{split}
\end{equation}
A discrete analogue of the Stokes-Dirac structure is the finite-dimensional Dirac structure constructed in the following theorem \cite{SeslijaCDC}.

\vspace{0.1cm}
\begin{thmm}\label{th:pdDirac}
Given linear spaces $\mathcal{F}_{p,q}^d$ and $\mathcal{E}_{p,q}^d$, and the bilinear form $\langle\!\langle, \rangle\!\rangle_d$. The linear subspace $\mathcal{D}_d\subset\mathcal{F}_{p,q}^d\times \mathcal{E}_{p,q}^d$ defined by
\begin{equation}\label{eq:Dir-prim-dual}
\begin{split}
& \mathcal{D}_{d}=\big\{ (\hat{f}_p, {f}_q, {f}_b,{e}_p,\hat{e}_q,\hat{e}_b)\in \mathcal{F}_{p,q}^d\times \mathcal{E}_{p,q}^d\big |\\
& \left(\begin{array}{c}\hat{f}_p \\ {f}_q \end{array} \right) =\left(\begin{array}{cc} 0  & (-1)^{r}  \mathbf{d}_\mathrm{i}^{n-q}\\  \mathbf{d}^{n-p} & 0\end{array} \right) \left(\begin{array}{c} {e}_p \\ \hat{e}_q\end{array} \right)+ (-1)^{r} \left(\begin{array}{c}  \mathbf{d}_\mathrm{b}^{n-q} \\ 0 \end{array}\right) \hat{e}_b\,,\\
& \begin{array}{c}~~~~~f_b \end{array}= ~(-1)^{p}\mathbf{tr}^{n-p}{e}_p \}\,,
 \end{split}
\end{equation}
with $r=pq+1$, is a Dirac structure with respect to the pairing $\langle\!\langle, \rangle\!\rangle_{d}$ .
\end{thmm}

Note that since $\mathbf{d}_\mathrm{i}^{n-q}=(-1)^q (\mathbf{d}^{n-p})^\textsc{t}$ and $\mathbf{d}_\mathrm{b}^{n-q}=(-1)^{n-p}(\mathbf{tr}^{n-p})^\textsc{t}$, the structure (\ref{eq:Dir-prim-dual}) is in fact a \emph{Poisson structure} on the state space $\Omega_d^p( \star_\mathrm{i}K)\times \Omega_d^q( K)$.


The other discrete analogue of the Stokes-Dirac structure is defined on the spaces
\begin{equation*}
\begin{split}
\tilde{\mathcal{F}}_{p,q}^d&=\Omega_d^p (K)\times \Omega_d^q( \star_\mathrm{i} K)\times \Omega_d^{n-p}(\partial (\star K))\\
\tilde{\mathcal{E}}_{p,q}^d&=\Omega_d^{n-p}( \star_\mathrm{i} K)\times \Omega_d^{n-q}( K)\times \Omega_d^{n-q}(\partial K)\,.
\end{split}
\end{equation*}
A natural discrete mirror of (\ref{eq-aj9}) in this case is a symmetric pairing defined by
\begin{equation*}
\begin{split}
\langle\!\langle (&\underbrace{{f}_p^1,\hat{f}_q^1,{\hat{f}}_b^1}_{\in \tilde{\mathcal{F}}_{p,q}^d},\underbrace{{\hat{e}}_p^1,{e}_q^1,{e}_b^1}_{\in \tilde{\mathcal{E}}_{p,q}^d}), ({f}_p^2,\hat{f}_q^2,\hat{f}_b^2,\hat{e}_p^2,{e}_q^2,{e}_b^2)\rangle\!\rangle_{\tilde{d}} \\
&= \langle \hat{e}_p^1\wedge {f}_p^2+ {e}_q^1\wedge {\hat{f}}_q^2+ {\hat{e}}_p^2\wedge {f}_p^1+ {e}_q^2\wedge \hat{f}_q^1,K \rangle\\
&~~\;+\langle {e}_b^1\wedge \hat{f}_b^2+ {e}_b^2\wedge \hat{f}_b^1,\partial K \rangle\,.
  \end{split}
\end{equation*}

\vspace{0.1cm}
\begin{thmm}\label{th:pdDiracN}
The linear space $\tilde{\mathcal{D}}_d$ defined by 
\begin{equation}\label{eq:Dir-prim-dual2}
\begin{split}
&\tilde{\mathcal{D}}_d=\big\{ ({f}_p, \hat{f}_q, {f}_b,{e}_p, {e}_q, {e}_b)\in\tilde{ \mathcal{F}}_{p,q}^d\times \tilde{\mathcal{E}}_{p,q}^d\big |\\
& \left(\!\!\begin{array}{c} {f}_p \\ {f}_q\end{array}\!\!\right)=\left(\!\!\begin{array}{cc}0 & (-1)^{pq+1}  \mathbf{d}^{n-q} \\  \mathbf{d}_\mathrm{i}^{n-p} & 0\end{array}\!\!\right)\left(\!\!\begin{array}{c} \hat{e}_p \\ {e}_q\end{array}\!\!\right)+\left(\!\!\begin{array}{c}  0\\ \mathbf{d}_\mathrm{b}^{n-p} \end{array}\!\!\right) \hat{f}_b\,,\\
&\begin{array}{c}\,\,~e_b ~\end{array}= ~(-1)^{p}\mathbf{tr}^{n-q}{e}_q\,
\big \}\,
 \end{split}
\end{equation}
is a Dirac structure with respect to the bilinear pairing $\langle\!\langle, \rangle\!\rangle_{\tilde{d}}$.
\end{thmm}

In the following section, the simplicial Dirac structures (\ref{eq:Dir-prim-dual}) and (\ref{eq:Dir-prim-dual2}) will be used as \emph{terminus a quo} for the geometric formulation of spatially discrete port-Hamiltonian systems.

\section{Port-Hamiltonian Systems on a Simplicial Complex}\label{Sec5}
Let a function $\mathcal{H}: \Omega_d^p(\star_\mathrm{i} K)\times \Omega_d^q(K)\rightarrow \mathbb{R}$ stand for the Hamiltonian 
$(\hat{\alpha}_p,\alpha_q)\mapsto \mathcal{H}(\hat{\alpha}_p,\alpha_q)$, with $\hat{\alpha}_p \in \Omega_d^p(\star_\mathrm{i} K)$ and $\alpha_q\in \Omega_d^q(K)$. A time derivative of $\mathcal{H}$ along an arbitrary trajectory $t\rightarrow (\hat{\alpha}_p(t),\alpha_q(t))\in  \Omega_d^p(\star_\mathrm{i} K)\times  \Omega_d^q( K)$, $t\in \mathbb{R}$, is
\begin{equation}
\begin{split}
\frac{\mathrm{d}}{\mathrm{d}t}\mathcal{H}(\hat{\alpha}_p,\alpha_q)=\langle \frac{\partial \mathcal{H}}{\partial \hat{\alpha}_p} \wedge \frac{\partial \hat{\alpha}_p}{\partial t} + \hat{\frac{\partial \mathcal{H}}{\partial \alpha_q}}\wedge \frac{\partial \alpha_q }{\partial t},K\rangle\,.
 \end{split}
\end{equation}
The relations between the simplicial-Dirac structure (\ref{eq:Dir-prim-dual}) and time derivatives of the variables are: $\hat{f}_p=-\frac{\partial \hat{\alpha}_p}{\partial t}$, $f_q=-\frac{\partial \alpha_q }{\partial t}$, while the efforts are: $e_p=\frac{\partial \mathcal{H}}{\partial \hat{\alpha}_p}$, $\hat{e}_q=\hat{\frac{\partial \mathcal{H} }{\partial \alpha_q}}$.

This allows us to define time-continuous port-Hamiltonian system on a simplicial complex $K$ (and its dual $\star K$) by
\begin{equation}\label{eq-38d}
\begin{split}
\!\!\!\!\!\!\!\!\!\!\!\!\!\!\!\!\!\!\!\!\!\!\left(\!\!\!\begin{array}{c}-\frac{\partial\hat{\alpha}_p}{\partial t} \\ -{\frac{\partial\alpha_q}{\partial t}}\end{array}\!\!\!\right)\!\!&=\!\!\left(\!\!\begin{array}{cc}0 & \!\!\!\!(-1)^{r}  \mathbf{d}_\mathrm{i}^{n-q}\\  \mathbf{d}^{n-p} & 0\end{array}\!\!\right)\!\!\left(\!\!\begin{array}{c} {\frac{\partial \mathcal{H}}{\partial \hat{\alpha}_p}} \\\hat{\frac{\partial \mathcal{H}}{\partial \alpha_q}}\end{array}\!\!\right)+\!(-1)^{r}\!\!\left(\!\!\begin{array}{c}  \mathbf{d}_\mathrm{b}^{n-q} \\ 0 \end{array}\!\!\right) \hat{e}_b\,,\!\!\!\!\\
\begin{array}{c}~~\,\,f_b ~~\end{array}&= ~(-1)^{p}\mathbf{tr}^{n-p}\frac{\partial \mathcal{H}}{\partial \hat{\alpha}_p}\,,
\end{split}
\end{equation}
where $r=pq+1$.

It immediately follows that $\frac{\mathrm{d}}{\mathrm{d}t}\mathcal{H}= \langle \hat{e}_b \wedge {f}_b, \partial K \rangle $, enunciating a fundamental property of the system: the increase in the energy on the domain $|K|$ is equal to the power supplied to the system through the boundary $\partial K$ and $\partial (\star K)$. The boundary efforts $\hat{e}_b$ are the boundary control input and $f_b$ are the outputs. 

\vspace{0.1cm}
\begin{rmk}\label{rmk:pass}
Introducing a linear negative feedback control as $\hat{e}_b=(-1)^{(n-p)(n-q)-1}*_\mathrm{b}f_b$, where $*_\mathrm{b}$ is the Hodge star on the boundary $\partial K$, leads to passivization of the lossless port-Hamiltonian system, i.e., $\frac{\mathrm{d}}{\mathrm{d}t}\mathcal{H}\leq-  \langle {f}_b \wedge *_\mathrm{b} {f}_b, \partial K \rangle \leq 0$. Furthermore, if the Hamiltonian is a $\mathcal{K}_\infty$ function with a strict minimum that is a stationary set for the system (\ref{eq-38d}), the equilibrium is asymptotically stable.
\end{rmk}

An alternative formulation of a spatially discrete port-Hamiltonian system is given in terms of the simplicial Dirac structure (\ref{eq:Dir-prim-dual2}). We start with the Hamiltonian function $
({\alpha}_p, \hat{\alpha}_q)\mapsto \mathcal{H}({\alpha}_p, \hat{\alpha}_q)$, where ${\alpha}_p \in \Omega_d^p(K)$ and $\hat{\alpha}_q\in \Omega_d^q(\star _\mathrm{i}K)$.
In a similar manner as in deriving (\ref{eq-38d}), we introduce the port-Hamiltonian system
\begin{equation}\label{eq-38dN}
\begin{split}
\!\!\!\!\!\!\!\left(\!\!\begin{array}{c}\!-\frac{\partial{\alpha}_p}{\partial t} \\ \!\!-{\frac{\partial \hat{\alpha}_q}{\partial t}}\!\!\!\end{array}\right)&\!=\!\left(\!\!\begin{array}{cc}0 & \!\!\!(-1)^{r}  \mathbf{d}^{n-q}\\  \mathbf{d}_\mathrm{i}^{n-p} & 0\end{array}\!\!\right)\!\left(\!\!\begin{array}{c} \hat{{\frac{\partial \mathcal{H}}{\partial{\alpha}_p}}} \\ {\frac{\partial \mathcal{H}}{\partial \hat{\alpha}_q}}\end{array}\!\!\right)\!+\!\left(\!\!\begin{array}{c} 0\\  \mathbf{d}_\mathrm{b}^{n-p}  \end{array}\!\!\right) \hat{f}_b\,,\\
\begin{array}{c}~~\,\,e_b ~~\end{array}&= ~(-1)^{p}\mathbf{tr}^{n-q}\frac{\partial \mathcal{H}}{\partial \hat{\alpha}_q}\,.
\end{split}
\end{equation}

In contrast to (\ref{eq-38d}), in the case of the formulation (\ref{eq-38dN}), the boundary flows $\hat{f}_b$ can be considered to be freely chosen, while the boundary efforts $e_b$ are determined by the dynamics. The free boundary variables are always defined on the boundary of the dual cell complex.


\section{Physical examples}
In this section we consider the discrete wave equation on a $2$-dimensional simplicial complex and the telegraph equations on a segment. 
\subsection{Two-dimensional wave equation}
Let us consider the simplicial Dirac structure underpining the discretized two-dimensional wave equation. The normalized wave equation is given by
$$\frac{\partial^2 \phi^c}{\partial t^2} - \Delta \phi^c=0\,,
$$
where $\phi^c$ is a smooth $0$-form on a compact surface $M\subset \mathbb{R}^2$ with a closed boundary, and $\Delta$ is the Laplace operator. (Throughout, the superscript $c$ designates the continuous quantities.) This equation, together with nonzero energy flow, can be formulated as a port-Hamiltonian system with boundary port variables \cite{Golo,Talasila}.

The energy variables of the discretized system are chosen as follows: the kinetic momentum is a dual $2$-form whose time derivative is set to be $\hat{f}_p$, the elastic strain is a primal $1$-form with time derivative corresponding to $f_q$, the coenergy variables are a primal $0$-form $e_p$ and a dual $1$-form $\hat{e}_q$. Such a formulation of the discrete wave equation is consonant with the simplicial Dirac structure (\ref{eq:Dir-prim-dual}) for the case when $p=n=2$ and $q=1$, and is given by\vspace{0.1cm}
\begin{equation*}
\begin{split}
 \left(\begin{array}{c}\hat{f}_p \\ {f}_q \end{array} \right) &=\left(\begin{array}{cc} 0  & -  \mathbf{d}_\mathrm{i}^1\\  \mathbf{d}^0 & 0\end{array} \right) \left(\begin{array}{c} {e}_p \\ \hat{e}_q\end{array} \right)- \left(\begin{array}{c}  \mathbf{d}_\mathrm{b}^1 \\ 0 \end{array}\right) \hat{e}_b\,,\\
 \begin{array}{c}f_b \end{array}&= ~\mathbf{tr}^0{e}_p\,.
 \end{split}
\end{equation*}

\begin{figure}
\centering
\vspace{-0.0cm}
      \hspace{0cm}\includegraphics[width=4.8cm]{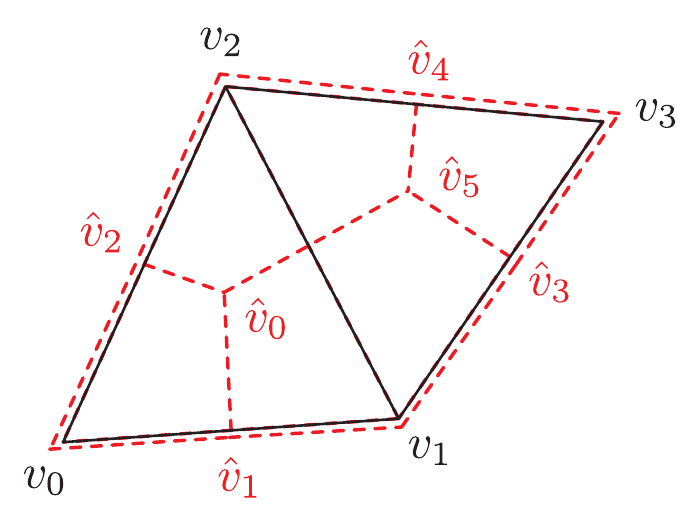}
      \vspace{-0.1cm}
  \caption{A simplicial complex $K$ consists of two triangles. The dual edges introduced by subdivision are shown dotted.}\label{fig:wave}
\end{figure}

\vspace{0.15cm}
\noindent\emph{Example.} Consider a simplicial complex pictorially given by Fig.~\ref{fig:wave}. The primal and dual $2$-faces have counterclockwise orientations. The matrix representation of the incidence operator $\partial_1$, from the primal edges to the primal vertices, is
\begin{equation*}
\left.\begin{array}{cccccc}
~ & [v_0,v_1] & [v_1,v_2] & [v_2,v_0] & [v_1,v_3] & [v_3,v_2] \\
v_0 & -1 & ~0 & ~0 & ~0 & ~0 \\
v_1 & ~1 & -1 & ~0 & -1 & ~0 \\
v_2 & ~0 & ~1 & -1 & ~0 & ~1 \\
v_3 & ~0 & ~0 & ~1 & ~1 & -1\end{array}\right.
\end{equation*}
while the discrete exterior derivative from the vertices to the edges is the transpose of the incidence operator, i.e., $\mathbf{d}^0=\partial_1^\textsc{t}$. The dual exterior derivative is $\mathbf{d}_\mathrm{i}^1=-\left(\mathbf{d}^0\right)^\textsc{t}$, while the matrix representation of the $\partial_{\mathrm{b},2}$ operator is
\begin{equation*}
\begin{array}{lcccc}
~ & \star_{\mathrm{i}} v_0 & \star_{\mathrm{i}} v_1 & \star_{\mathrm{i}} v_2 & \star_{\mathrm{i}} v_3\\
\left[\hat{v}_2,\hat{v}_1\right] & 1 & 0 & 0 & 0 \\
\left[\hat{v}_1,\hat{v}_3\right] & 0 & 1 & 0 & 0  \\
\left[\hat{v}_3,\hat{v}_4\right] & 0 & 0 & 0 & 1 \\
\left[\hat{v}_4,\hat{v}_2\right] & 0 & 0 & 1 & 0\end{array}
\end{equation*}
The trace operator is $\mathbf{tr}^0=  ( \mathbf{d}_\mathrm{b}^1 )^\textsc{t}=\partial_{\mathrm{b},2}$.

It is trivial to show
\begin{equation}
\begin{split}
\langle &\mathbf{d}^0 e_p \wedge \hat{e}_ q ,K \rangle + \langle  e_p \wedge ( \mathbf{d}_\mathrm{i}^1 \hat{e}_ q +\mathbf{d}_\mathrm{b}^1 \hat{e}_ b ),K \rangle \\
&= \hat{e}_b[\hat{v}_2,\hat{v}_1] f_b(v_0) + \hat{e}_b[\hat{v}_1,\hat{v}_3] f_b(v_1)\\
&~+\hat{e}_b[\hat{v}_3,\hat{v}_4] f_b(v_3)+\hat{e}_b[\hat{v}_4,\hat{v}_2] f_b(v_2)\,,
 \end{split}
\end{equation}
what confirms that the boundary terms genuinely live on the boundary of $|K|$.

\subsection{Telegraph equations}
We consider an ideal lossless transmission line on a $1$-dimensional simplicial complex. The energy variables are the charge density ${q}\in \Omega_d^1(K)$, and the flux density $\hat{\phi}\in\Omega_d^1(\star  K)$, hence $p=q=1$. The Hamiltonian representing the total energy stored in the transmission line with distributed capacitance $C$ and distributed inductance $\hat{L}$ is\vspace{-0.1cm}
\begin{equation}
\begin{split}
\mathcal{H}=  \langle \frac{1}{2C} {q} \wedge *_1 {q}  +  \frac{1}{2\hat{L}}  \hat{\phi} \wedge *_0^{-1} \hat{\phi}  ,K\rangle \,,
 \end{split}
\end{equation}
where $*_0$ and $*_1$ are the discrete diagonal Hodge operators that relate the appropriate cochains according to the following schematic diagram
\begin{equation*}
\left.\begin{array}{ccccc} \!\!\! \Omega_d^0(\partial K)\!\!\! & \!\!\! \xleftarrow{~\mathbf{tr}^0~} \!\!\! & \!\!\! \Omega_d^0(K)\!\!\! & \!\xrightarrow{~\mathbf{d}^0~}\!\! &  \Omega_d^1(K)\!\!\! \\ 
\!\!\!\downarrow\! {*_\mathrm{b}} \!\!\!& \!\!\! \!\!\! &\!\!\!  \downarrow\!{*_0} \!\!\!&\!\!\!\!\!\!  & \!\!\! \downarrow\!{*_1} \!\!\! \\ 
\!\!\! \Omega_d^0(\partial(\star  K))\! & \xrightarrow{~\mathbf{d}_\mathrm{b}^0~} &  \Omega_d^1(\star_\mathrm{i} K)  & \!\xleftarrow{{~\mathbf{d}}_\mathrm{i}^0~}\! & ~ \Omega_d^0(\star_\mathrm{i} K)\,.\!\!\! \end{array}\right.\,
\end{equation*}

The co-energy variables are: $\hat{e}_p=\hat{\frac{\partial \mathcal{H}}{\partial{q}}}=*\frac{{q}}{C}=\hat{V}$ representing voltages and ${e}_q=\frac{\partial \mathcal{H}}{\partial \hat{\phi }}=*\frac{\hat{\phi}}{\hat{L}}=I$ currents. Selecting ${f}_p=-\frac{\partial{q}}{\partial  t}$ and $\hat{f}_q=-\frac{\partial \hat{\phi}}{\partial  t}$ leads to the port-Hamiltonian formulation of the telegraph equations
\begin{equation}\label{eq:Dir-prim-dualMET}
\begin{split}
\left(\begin{array}{c} - \frac{\partial {q} }{\partial t} \\ -\frac{\partial \hat{\phi} }{\partial t}\end{array}\right)&=\left(\begin{array}{cc} 0 & \mathbf{d}^0 \\  \mathbf{d}_\mathrm{i}^0 & 0\end{array}\right) \left(\begin{array}{c} *_1\frac{q}{C} \\  *_0^{-1}\frac{ \hat{\phi} }{\hat{L}}  \end{array}\right)+\left(\begin{array}{c} 0 \\  \mathbf{d}_\mathrm{b}^0  \end{array}\right) \hat{f}_b\,\\
\begin{array}{c}~~\,\,e_b ~~\end{array} &= ~ - \mathbf{tr}^0\left(*_0^{-1} \frac{\hat{\phi}}{\hat{L}}\right) \,,
 \end{split}
\end{equation}where $\hat{f}_b$ are the input voltages and $e_b$ are the output currents.

\begin{figure}
\centering
    \includegraphics[width=7cm]{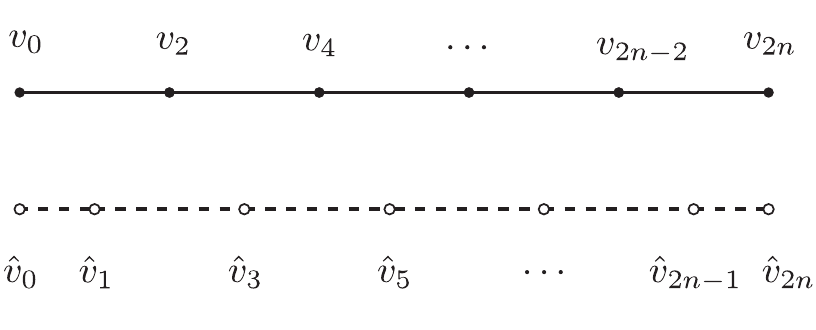}\\\vspace{-0.4cm}
  \caption{The primal $1$-dimensional simplicaial complex $K$. By construction, the nodes $\hat{v}_0$ and $\hat{v}_{2n}$ are added to the boundary to insure that $\partial (\star K)=\star (\partial K)$. }\label{fig:RDcomp1}
\end{figure}

In the case we wanted to have the electrical currents as the inputs, the charge and the flux densities would be defined on the dual mesh and the primal mesh, respectively. Instead of the port-Hamiltonian system in the form (\ref{eq:Dir-prim-dualMET}), the discretized telegraph equations would be in the form (\ref{eq-38d}). The charge density is defined on the dual cell complex as $\hat{q}\in \Omega_d^1(\star_\mathrm{i} K)$ and the discrete flux density is $\phi \in \Omega_d^1(K)$. The finite-dimensional port-Hamiltonian system is of the form
\begin{equation}\label{eq:Dir-prim-dualMETcurrent}
\begin{split}
\left(\begin{array}{c}- \frac{\partial \hat{q} }{\partial t} \\ -\frac{\partial \phi }{\partial t}\end{array}\right)&=\left(\begin{array}{cc} 0 & \mathbf{d}_\mathrm{i}^0 \\  \mathbf{d}^0 & 0\end{array}\right) \left(\begin{array}{c} *_0^{-1} \frac{\hat{q}}{\hat{C}} \\  *_1\frac{ \phi }{L}  \end{array}\right) +\left(\begin{array}{c}   \mathbf{d}_\mathrm{b}^0 \\ 0 \end{array}\right) \hat{e}_b\,\\
\begin{array}{c}~~\,\,f_b ~~\end{array} &= ~ - \mathbf{tr}^0 \left(*_0^{-1} \frac{\hat{q}}{\hat{C}} \right) \,,
 \end{split}
\end{equation}
where $\hat{e}_b$ are the input currents and $f_b$ are the output voltages.

\vspace{0.15cm}
The exterior derivative $\mathbf{d}^0: \Omega_d^0(K)\rightarrow \Omega_d^1(K)$ is the transpose of the incidence matrix of the primal mesh. The discrete derivative $\mathbf{d}_\mathrm{i}^0: \Omega_d^0(\star_\mathrm{i} K)\rightarrow \Omega_d^1(\star_\mathrm{i} K)$ in the matrix notation is the incidence matrix of the primal mesh. Thus, we have
\begin{equation}
\begin{split}
(\mathbf{d}_\mathrm{i}^0)^\textsc{t}=\mathbf{d}^0&=\left(\begin{array}{rrrrrr}-1 & 1 & 0 & \cdots & 0 & 0 \\0 & -1 & 1 & \cdots & 0 & 0 \\ \, & \, & \, & \ddots & \, & \, \\0 & 0 & 0 & \cdots & -1 & 1\end{array}\right)\,,\\
\mathbf{tr}^0=(\mathbf{d}_\mathrm{b}^0)^\textsc{t}&=\left(\begin{array}{rrrrrr}-1 & 0 & 0 & \cdots & 0 & 0 \\0 & 0 & 0 & \cdots & 0 & 1\end{array}\right)\,.
\end{split}
\end{equation}


\vspace{0.1cm}
\begin{rmk}The discrete analogue of the Stokes-Dirac structure obtained in \cite{Golo} is a finite-dimensional Dirac structure, but not a Poisson structure. The implication of this on the physical realization is that the transmission line in the finite-dimensional case is not only composed of inductors and capacitors but also of transformers.\end{rmk}

\begin{figure}
\centering
    \includegraphics[width=7.8cm]{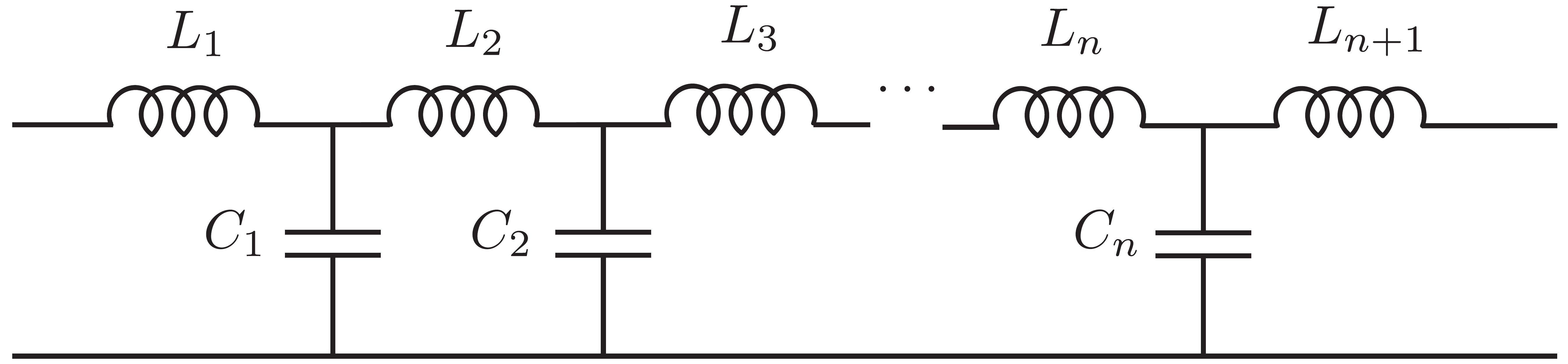}\\\vspace{-0.0cm}
  \caption{The finite-dimensional approximation of the lossless transmission line when the inputs are voltages and the outputs are currents. The inductances $L_1,\ldots, L_{n+1}$ are the values that the discrete distributed inductance $\hat{L}$ takes on the simplices $[\hat{v}_0,\hat{v}_1],\ldots, [\hat{v}_{2n-1},\hat{v}_{2n}]$; the capacitances $C_1,\ldots C_n$ are the values $C$ takes on $[v_0,v_2],\ldots, [v_{2n-2},v_{2n}]$.}\label{fig:LC1}
  \end{figure}

\begin{figure}
\centering
    \includegraphics[width=7.8cm]{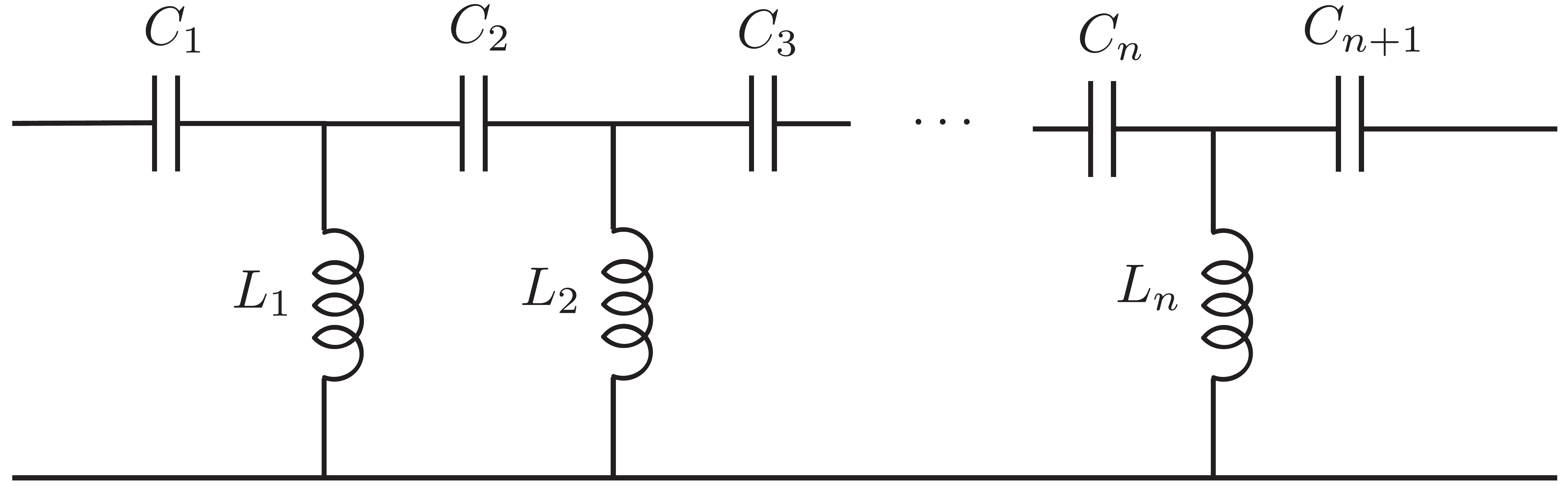}\\\vspace{-0.0cm}
  \caption{The finite-dimensional approximation of the lossless transmission line when the inputs are currents and the outputs are voltages. The inductances are: $L_1=\int_{[v_0,v_2]}L^c=L([v_0,v_2])$, $L_2=\int_{[v_2,v_4]}L^c=L([v_2,v_4])$, $\ldots$, $L_n=\int_{[v_{2n-2},v_{2n}]} L^c =L([v_{2n-2},v_{2n}])$; the values of capacitors are: $C_1=\int_{[\hat{v}_0,\hat{v}_1]}C^c=\hat{C}([\hat{v}_0,\hat{v}_1])$, $C_2=\int_{[\hat{v}_1,\hat{v}_3]}C^c=\hat{C}([\hat{v}_1,\hat{v}_3])$, $C_3=\int_{[\hat{v}_3,\hat{v}_5]}C^c=\hat{C}([\hat{v}_3,\hat{v}_5])$, $\ldots$, $C_{n+1}=\int_{[\hat{v}_{2n-1},\hat{v}_{2n}]}C^c=\hat{C}([\hat{v}_{2n-1},\hat{v}_{2n}])$.}\label{fig:LC2}
  \end{figure}

The physical realizations of the port-Hamiltonian systems (\ref{eq:Dir-prim-dualMET}) and (\ref{eq:Dir-prim-dualMETcurrent}) are given on  Fig.~\ref{fig:LC1} and Fig.~\ref{fig:LC2}, respectively. Stabilization of either of those systems is easily achieved by terminating boundary ports with resistive elements, what is a practical application of the passivization explained in Remark~\ref{rmk:pass}.

The accuracy of the proposed method is $1/n$ (see \cite{SeslijaArXiv}).


\end{document}